\documentclass{article}

\usepackage[applemac]{inputenc} % Spécifique au LaTeX sous Macintosh
\usepackage[dvips]{graphicx}
\usepackage{amsmath} % Extra environments for multiline displayed eq
\usepackage{amssymb} % Defines symbol names 

\usepackage{colortbl}
\usepackage{soul} % The package soul is needed for \highlighttext to work.

\usepackage{pstcol}
\usepackage{pst-node}
\usepackage{pstricks}
\usepackage{fancybox}
\usepackage{calc}
\usepackage{truncate}
\usepackage{comment}

\usepackage{graphics} % nécessaire avec le package listings
\usepackage{listings}    % pour lister des programmes

\makeatletter
\def\@maketitle{%
  \vbox to 2.3in{%
    \hsize\textwidth
    \linewidth\hsize
    \vspace*{1.5cm}
    \centering
    {\bfseries\huge \@title \par}
    \vskip 2em
    {\large \begin{tabular}[t]{c}\@author \end{tabular}\par}
    \vfill}    \vspace*{1.0cm}
}
\renewcommand\section{\@startsection {section}{1}{\z@}%
     {.7\baselineskip plus\baselineskip}{.5\baselineskip}
                                   {\normalfont\Large\bfseries}}
\renewcommand\section{\@startsection {section}{1}{\z@}%
      {.5\baselineskip\@plus.7\baselineskip}{.3\baselineskip}%
                                   {\normalfont\Large\bfseries}}
\renewcommand\subsection{\@startsection{subsection}{2}{\z@}%
       {.5\baselineskip\@plus.7\baselineskip}{.3\baselineskip}%
                                   {\normalfont\large\bfseries}}
\renewcommand\subsubsection{\@startsection{subsubsection}{3}{\z@}%
      {.5\baselineskip\@plus.7\baselineskip}{.3\baselineskip}%
                                     {\normalfont\normalsize\bfseries}}
\makeatother

\renewenvironment{abstract}%
  {\normalfont
    \list{}{\labelwidth0pt
      \leftmargin0pt \rightmargin\leftmargin
      \listparindent\parindent \itemindent0pt
      \parsep0pt
      
    }%
    \item[\hskip\labelsep\bfseries\abstractname\enspace --] \itshape%
}{%
  \endlist}

\newcommand{\keywordsname}{Keywords}
\newenvironment{keywords}%
  {\normalfont
    \list{}{\labelwidth0pt
      \leftmargin0pt \rightmargin\leftmargin
      \listparindent\parindent \itemindent0pt
      \parsep0pt
      }%
    \item[\hskip\labelsep\bfseries\keywordsname:]}{\endlist}

\begin{document}

\pagestyle{myheadings}

\title{
On the generation of hyper-powersets\\
for the DSmT}

\author{Jean Dezert\\
ONERA\\
29 Av. de la  Division Leclerc \\
92320 Ch\^{a}tillon, France.\\
Jean.Dezert@onera.fr\\
\and
Florentin Smarandache\\
Department of Mathematics\\
University of New Mexico\\
Gallup, NM 87301, U.S.A.\\
smarand@unm.edu}
\date{}

\maketitle
%\vspace{2cm}

\begin{abstract}
The recent theory of plausible and paradoxical reasoning (DSmT) developed by the authors appears to be a nice promising theoretical tools to solve many information fusion problems where the Shafer's model cannot be used due to the intrinsic paradoxical nature of the elements of the frame of discernment and where a strong internal conflict between sources arises. The main idea of DSmT is to work on the hyper-powerset of the frame of discernment of the problem under consideration. Although the definition of hyper-powerset is well established, the major difficulty in practice is to generate such hyper-powersets in order to implement DSmT fusion rule on computers. We present in this paper a simple algorithm for generating hyper-powersets and discuss the limitations of our actual computers to generate such hyper-powersets when the dimension of the problem increases.
\end{abstract}

\begin{keywords}
Dezert-Smarandache theory (DSmT), hyper-powersets, monotone boolean functions, Dedekind problem, plausible and paradoxical reasoning, data fusion.
\end{keywords}

%--------------------------
\section{Introduction}
%--------------------------

The Dezert-Smarandache theory (DSmT for short) of plausible and paradoxical reasoning \cite{Dezert_2002b,Dezert_2003,Smarandache_2002} is a generalization of the classical Dempster-Shafer theory (DST) \cite{Shafer_1976} which allows to formally combine any types of sources of information (rational, uncertain or paradoxical). The DSmT is able to solve complex fusion problems where the DST usually fails, specially when conflicts (paradoxes) between sources become large and when the refinement of the frame of discernment $\Theta$ is inaccessible because of the vague, relative and imprecise nature of elements of $\Theta$ (see \cite{Dezert_2003} for justification and examples). The foundation of DSmT is based on the definition of the hyper-powerset $D^\Theta$ of a general frame of discernment $\Theta$.  $\Theta$ must be considered as a set $\{\theta_{1},\ldots,\theta_{n}\}$ of $n$  elements considered as exhaustive which cannot be precisely defined and separated so that no refinement of $\Theta$ in a new larger set $\Theta_{ref}$ of disjoint elementary hypotheses is possible in contrast with the classical Dempster-Shafer Theory (DST).  The DSmT deals directly with paradoxical/conflicting sources of information into this new formalism and proposes a new and very simple (associative and commutative) rule of combination for conflicting sources of informations (corpus/bodies of evidence). Some interesting results based on DSmT approach can be found in \cite{Tchamova_2003,Corgne_2003}.
Before going deeper into the DSmT it is necessary to briefly present first the foundations of the DST and DSmT for a better understanding of the important differences between these two theories.

%*****************************************
\section{Short presentation of the DST}
%*****************************************

Let $\Theta=\{\theta_1,\theta_2,\ldots,\theta_n\}$ be the frame of discernment of the problem under consideration having $n$ {\it{exhaustive}} and {\it{exclusive elementary}} hypothesis $\theta_i$. This corresponds to the Shafer's model of the problem.  Such model assumes that an ultime refinement of the problem is possible so that $\theta_i$ are well precisely defined/identified in such a way  that we are sure that they are exclusive and exhaustive. From this Shafer's model,  a basic belief assignment (bba) $m (.): 2^\Theta \rightarrow  [0, 1]$  associated to a given body of evidence $\mathcal{B}$ (also called sometimes corpus of evidence) is defined by
\begin{equation}
m(\emptyset)=0  \qquad \text{and}\qquad     \sum_{A\in 2^\Theta} m(A) = 1                            
%\end{equation}
%\begin{equation}
%\sum_{A\in 2^\Theta} m(A) = 1
\end{equation}
\noindent
where $2^\Theta$ is called the {\it{powerset}} of $\Theta$, i.e. the set of all subsets of $\Theta$. From any bba, one defines the belief and plausibility functions of $A\subseteq\Theta$ as
\begin{equation}
\text{Bel}(A) = \sum_{B\in 2^\Theta, B\subseteq A} m(B)
\label{Belg}
\end{equation}
\begin{equation}
\text{Pl}(A) = \sum_{B\in 2^\Theta, B\cap A\neq\emptyset} m(B)=1- \text{Bel}(\bar{A})
\label{Plg}
\end{equation}

Now let $\text{Bel}_1(.)$ and $\text{Bel}_2(.)$ be two belief functions over the same frame of discernment $\Theta$ and their corresponding bba $m_1(.)$ and $m_2(.)$ provided by two distinct bodies of evidence $\mathcal{B}_1$ and $\mathcal{B}_2$. Then the combined global belief function denoted $\text{Bel}(.)= \text{Bel}_1(.)\oplus \text{Bel}_2(.)$ is obtained by combining the information granules $m_1(.)$ and $m_2(.)$ through the following Dempster's rule of combination $[m_{1}\oplus m_{2}](\emptyset)=0$ and $\forall B\neq\emptyset \in 2^\Theta$,
 \begin{equation}
[m_{1}\oplus m_{2}](B) = 
\frac{\sum_{X\cap Y=B}m_{1}(X)m_{2}(Y)}{1-\sum_{X\cap Y=\emptyset} m_{1}(X) m_{2}(Y)} 
\label{eq:DSR}
 \end{equation}
 
The notation $\sum_{X\cap Y=B}$ represents the sum over all $X, Y \in 2^\Theta$ such that $X\cap Y=B$. 
The orthogonal sum $m (.)\triangleq [m_{1}\oplus m_{2}](.)$ is considered as a basic belief assignment if and only if the denominator in equation \eqref{eq:DSR} is non-zero. The term $k_{12}\triangleq \sum_{X\cap Y=\emptyset} m_{1}(X) m_{2}(Y)$ is called degree of conflict between the sources $\mathcal{B}_1$ and $\mathcal{B}_2$. When $k_{12}=1$,  the orthogonal sum $m (.)$ does not exist and the bodies of evidences $\mathcal{B}_1$ and $\mathcal{B}_2$ are said to be in {\it{full contradiction}}. Such a case can arise when there exists $A \subset \Theta$ such that $\text{Bel}_1(A) =1$ and $\text{Bel}_2(\bar{A}) = 1$. Same kind of trouble can occur also with the {\it{Optimal Bayesian Fusion Rule}} (OBFR) \cite{Dezert_2001a,Dezert_2001b}.\\

The DST is attractive for the {\it{Data Fusion community}} because it gives a nice mathematical model for ignorance and it includes the Bayesian theory as a special case \cite{Shafer_1976} (p. 4). Although very appealing, the DST presents some weaknesses and limitations because of its model itself, the theoretical justification of the Dempster's rule of combination but also because of our confidence to trust the result of Dempster's rule of combination when the conflit becomes important between sources ($k_{12} \nearrow 1$).\\

The Dempster's rule of combination has however been {\it{a posteriori}} justified  by the Smet's axiomatic of the Transferable Belief Model (TBM) in \cite{Smets_1994}. But we must also emphasize here that an infinite number of possible rules of combinations can be built from the Shafer's model  following ideas initially proposed in \cite{Lefevre_2002} and corrected here as follows:
\begin{itemize}
\item one first has to compute $m(\emptyset)$ by
$$ m(\emptyset) \triangleq \sum_{A\cap B =\emptyset}m_1(A)m_2(B) $$
\item then one redistributes $m(\emptyset)$ on all $(A\neq\emptyset)\subseteq \Theta$ with some given coefficients $w_m(A)\in[0,1]$ such that $\sum_{A\subseteq \Theta} w_m(A)=1$ according to
\begin{equation}
\begin{cases}
w_m(\emptyset)m(\emptyset) \rightarrow m(\emptyset)\\
m(A) + w_m(A)m(\emptyset) \rightarrow m(A), \forall A\neq\emptyset
\end{cases}
\label{eq:CEV}
\end{equation}
\end{itemize}
The particular choice of the set of coefficients $w_m(.)$ provides a particular rule of combination.
Actually there exists an infinite number of possible rules of combination. Some rules can be better justified than others depending on their ability or not to preserve associativity and commutativity properties of the combination. It can be easily shown in  \cite{Lefevre_2002}  that such general procedure provides all existing rules developed in the literature from the Shafer's model as alternative to the primeval Dempster's rule of combination depending on the choice of coefficients $w(A)$. As example the Dempster's rule of combination can be obtained from \eqref{eq:CEV} by choosing $w_m(\emptyset)=0$ and $ w_m(A)=m(A)/(1-m(\emptyset))$ for all $A\neq\emptyset$. The Yager's rule of combination is obtained by choosing $w_m(\Theta)=1$ while the "Smets' rule of combination" is obtained by choosing $w_m(\emptyset)=1$ and thus accepting the possibility to deal with bba such that $m(\emptyset)>0$.

%*****************************************
\section{Foundations of the  DSmT}
%*****************************************

The development of the Dezert-Smarandache theory of plausible and paradoxical reasoning (called DSmT for short) comes from the necessity to overcome the two following inherent limitations of the DST which are closely related with the acceptance of the third middle excluded principle, i.e.
\begin{enumerate}
\item[(C1)] - the DST considers a discrete and finite frame of discernment 
$\Theta$ based on a set of exhaustive and exclusive elementary elements $\theta_i$.
\item[(C2)] - the bodies of evidence are assumed independent and provide their own belief function on the powerset $2^\Theta$ but with {\it{same interpretation}} for $\Theta$. 
\end{enumerate}

The foundation of the DSmT is based on the refutation of the principle of the third excluded middle for a wide class of fusion problems due to the nature of the elements of $\Theta$. By accepting the third middle, we can easily handle the possibility to deal directly with paradoxes (partial vague overlapping elements/concepts) of the frame of discernment. This is the DSm model. In other words, we include the third exclude directly into the formalism to develop the DSmT and relax the (C1) and (C2) constraints of the Shafer's model. By doing this, a wider class of fusion problem can be attacked by the DSmT.
The relaxation of the constraint (C1) can be justified since, in many problems (see example in \cite{Dezert_2003}), the elements of $\Theta$ generally correspond only to imprecise/vague notions and concepts so that no refinement of $\Theta$ satisfying the first constraint is actually possible (specially if natural language is used to describe elements of $\Theta$).\\

The DSmT refutes also the excessive requirement imposed by (C2) in the Shafer's model, since it seems clear to us that, the {\sl{same}} frame $\Theta$ may be interpreted differently by the distinct bodies of evidence (experts). Some subjectivity on the information provided by a source of information is almost unavoidable, otherwise this would assume, as within the DST, that all bodies of evidence have an objective/universal (possibly uncertain) interpretation or measure of the phenomena under consideration which unfortunately rarely occurs in reality, but when bba are based on some {\it{objective probabilities}} transformations. But in this last case, probability theory can handle properly the information; and the DST, as well as the DSmT, becomes useless. If we now get out of the probabilistic background argumentation, we claim that in most of cases, the sources of evidence provide their beliefs about some hypotheses only with respect to their own worlds of knowledge and experience without reference to the (inaccessible) absolute truth of the space of possibilities.
The DSmT  includes the possibility to deal with evidences arising from different sources of information which don't have access to absolute interpretation of the elements $\Theta$ under consideration. The DSmT can be interpreted as a general and direct extension of Bayesian theory and the Dempster-Shafer theory in the 
following sense. Let $\Theta=\{\theta_{1},\theta_{2}\}$ be the simpliest frame of 
discernment involving only two elementary hypotheses (with no  additional assumptions on $\theta_{1}$ and $\theta_{2}$), then 
\begin{itemize}
\item the probability theory deals with basic probability assignments (bpa)
$m(.)\in [0,1]$ such that $m(\theta_{1})+m(\theta_{2})=1$
\item the DST deals with bba $m(.)\in [0,1]$ such that 
$m(\theta_{1})+m(\theta_{2})+m(\theta_{1}\cup\theta_{2})=1$
\item the DSmT theory deals with
generalized bba $m(.)\in [0,1]$ such that 
$m(\theta_{1})+m(\theta_{2})+m(\theta_{1}\cup\theta_{2})+m(\theta_{1}\cap\theta_{2})=1$
\end{itemize}

%*****************************************************
\subsection{Notion of hyper-powerset $D^\Theta$}
%*****************************************************

One of the cornerstones of the DSmT is  the notion of hyper-powerset which is now presented. Let $\Theta=\{\theta_{1},\ldots,\theta_{n}\}$ be a set of $n$ 
elements which cannot be precisely defined and separated so that no 
refinement of $\Theta$ in a new larger set $\Theta_{ref}$ of disjoint elementary 
hypotheses is possible (we abandon here the Shafer's model). The {\sl{hyper-powerset}} $D^\Theta$ is defined as the set of all composite propositions built from elements of $\Theta$ with $\cup$ and $\cap$ ($\Theta$ generates $D^\Theta$ under operators $\cup$ and $\cap$)
operators such that 
\begin{enumerate}
\item $\emptyset, \theta_1,\ldots, \theta_n \in D^\Theta$.
\item  If $A,B \in D^\Theta$, then $A\cap B\in D^\Theta$ and $A\cup B\in D^\Theta$.
\item No other elements belong to $D^\Theta$, except those obtained by using rules 1 or 2.
\end{enumerate}
The dual (obtained by switching $\cup$ and $\cap$ in expressions) of $D^\Theta$ is itself.  There are elements in $D^\Theta$ which are self-dual (dual to themselves), for example $\alpha_8$ for the case when $n=3$ in the example below.
The cardinality of $D^\Theta$ is majored by 
$2^{2^n}$ when $\text{Card}(\Theta)=\vert\Theta\vert=n$. The generation 
of hyper-power set $D^\Theta$ is closely related with the famous Dedekind's problem \cite{Dedekind_1897,Comtet_1974} on enumerating the 
set of monotone Boolean functions as it will be presented in the sequel with  the generation of the elements of $D^\Theta$.\\

\noindent{\it{Example of the first hyper-powersets $D^\Theta$}}

In the degenerate case ($n=0)$ where $\Theta=\{  \}$, one has $D^\Theta=\{\alpha_0\triangleq\emptyset\}$ and $\vert D^\Theta\vert = 1$. When $\Theta=\{\theta_{1}\}$, one has $D^\Theta=\{\alpha_0\triangleq\emptyset,\alpha_1\triangleq\theta_1 \}$ and $\vert D^\Theta\vert = 2$. When $\Theta=\{\theta_{1},\theta_{2}\}$, one has
$D^\Theta=\{\alpha_0,\alpha_1,\ldots,\alpha_{4} \}$ and $\vert D^\Theta\vert = 5$ with
$\alpha_0\triangleq\emptyset$, $\alpha_1\triangleq\theta_1\cap\theta_2$, $\alpha_2\triangleq\theta_1$, $\alpha_3\triangleq\theta_2 $ and $\alpha_4\triangleq\theta_1\cup\theta_2 $. When $\Theta=\{\theta_{1},\theta_{2},\theta_{3}\}$,  one has $D^\Theta=\{\alpha_0,\alpha_1,\ldots,\alpha_{18} \}$ and $\vert D^\Theta\vert = 19$ (see \cite{Dezert_2003} for details) with
\begin{equation*}
\begin{array}{l}
\alpha_i     \\
\hline
\alpha_0\triangleq\emptyset                                                 \\
\alpha_1\triangleq\theta_1\cap\theta_2\cap\theta_3       \\
\alpha_2\triangleq\theta_1\cap\theta_2                                \\
\alpha_3\triangleq\theta_1\cap\theta_3                                \\
\alpha_4\triangleq\theta_2\cap\theta_3                                \\
\alpha_5\triangleq(\theta_1\cup\theta_2)\cap\theta_3       \\
\alpha_6\triangleq(\theta_1\cup\theta_3)\cap\theta_2       \\
\alpha_7\triangleq(\theta_2\cup\theta_3)\cap\theta_1       \\
\alpha_8\triangleq[(\theta_1\cap\theta_2)\cup\theta_3] \cap(\theta_1\cup\theta_2)      \\
\alpha_9\triangleq\theta_1                                                      \\
\alpha_{10}\triangleq\theta_2                                                     \\
\alpha_{11}\triangleq\theta_3                                                      \\
\alpha_{12}\triangleq(\theta_1\cap\theta_2)\cup\theta_3      \\
\alpha_{13}\triangleq(\theta_1\cap\theta_3)\cup\theta_2      \\
\alpha_{14}\triangleq(\theta_2\cap\theta_3)\cup\theta_1      \\
\alpha_{15}\triangleq(\theta_1\cup\theta_2)                             \\
\alpha_{16}\triangleq(\theta_1\cup\theta_3)                             \\
\alpha_{17}\triangleq(\theta_2\cup\theta_3)                             \\
\alpha_{18}\triangleq(\theta_1\cup\theta_2\cup\theta_3)       \\
\end{array}
\end{equation*}

Note that the {\it{classical}} complementary $\bar{A}$ of any proposition $A$ (except for $\emptyset$ and $\Theta$), is not involved within DSmT because of the refutation of the third excluded middle.  $\vert D^\Theta\vert$ for $n\geq1$ follows the sequence of Dedekind's numbers\footnote{Actually this sequence corresponds to the sequence of Dedekind minus one since we don't count the last degenerate isotone function $f_{2^{2^n}-1}(.)$ as element of $D^\Theta$ (see section 4.2)} 1,2,5,19,167,7580,7828353,... \cite{Sloane_2003}.\\

From a general frame of discernment $\Theta$, we define a map $m(.): 
D^\Theta \rightarrow [0,1]$ associated to a given body of evidence $\mathcal{B}$ 
which can support paradoxical information, as follows
\begin{equation*}
m(\emptyset)=0 \qquad \text{and}\qquad \sum_{A\in D^\Theta} m(A) = 1 
\end{equation*}
The quantity $m(A)$ is called $A$'s {\sl{generalized basic belief 
assignment}} (gbba) or the generalized basic belief mass for $A$.
The belief and plausibility functions are defined in almost the same manner as within the DST, i.e.
\begin{equation*}
\text{Bel}(A) = \sum_{\overset{B\subseteq A}{B\in D^\Theta}} m(B)
\qquad\text{and}\qquad\text{Pl}(A) = \sum_{\overset{B\cap A\neq\emptyset}{B\in D^\Theta}} m(B)
\end{equation*}

These definitions are compatible with the DST definitions when the sources of 
information become uncertain but rational (they do not support paradoxical 
information). We still have $\forall A\in D^\Theta, \text{Bel}(A)\leq \text{Pl}(A)$.

%***************************************
\subsection{The DSm rule of combination}
%***************************************

The  DSm rule of combination $m(.)\triangleq [m_{1}\oplus m_{2}](.)$ of two distinct (but potentially paradoxical) sources of evidences $\mathcal{B}_{1}$ and  $\mathcal{B}_{2}$ over the same 
general frame of discernment $\Theta$ with belief functions $\text{Bel}_{1}(.)$ and 
 $\text{Bel}_{2}(.)$ associated with general information granules $m_{1}(.)$ and $m_{2}(.)$ is  given by $\forall C\in D^\Theta$,
 \begin{equation}
m(C) = 
 \sum_{A,B\in D^\Theta, A\cap B=C}m_{1}(A)m_{2}(B)
 \label{JDZT}
 \end{equation}
Since $D^\Theta$ is closed under $\cup$ and $\cap$ operators, this new rule 
of combination guarantees that $m(.): D^\Theta \rightarrow [0,1]$ is a proper general information granule. This rule of combination is commutative and associative 
and can always be used for the fusion of paradoxical or rational sources of 
information (bodies of evidence). It is important to note that any fusion of sources of information 
generates either uncertainties, paradoxes or {\sl{more generally both}}. This is intrinsic to the 
general fusion process itself. The theoretical justification of the DSm rule can be found in \cite{Dezert_2003}.\\

This DSm rule of combination seems at the first glance very expensive in terms of computations and memory size due to the huge number of elements in $D^\Theta$. This is only true if the cores (the set of focal elements of gbba) $\mathcal{K}_1(m_1)$ and $\mathcal{K}_2(m_2)$ coincide with $D^\Theta$; in other words when $m_1(A)>0$ and $m_2(A)>0$ for all $A\neq\emptyset\in D^\Theta$. Fortunately, it is important to note here that in most of practical applications the sizes of $\mathcal{K}_1(m_1)$ and $\mathcal{K}_2(m_2)$ are much smaller than $\vert D^\Theta\vert$ because bodies of evidence generally allocate their basic belief assignments only over a subset of hyper-powerset. This makes things easier for the implementation of the DSm rule \eqref{JDZT}. 
The DSm rule is actually very easy to implement. It suffices for each focal element of $\mathcal{K}_1(m_1)$ to multiply it with the focal elements of $\mathcal{K}_2(m_2)$ and then to  pool all combinations which are logically equivalent under the  algebra of sets according to the following scheme
\begin{figure}[hbtp]
\centering
 \includegraphics[height=4cm]{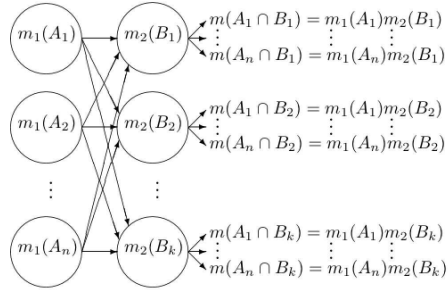}
\caption{Representation of the DSm rule}
 \label{fig:1}
 \end{figure}

The figure above represents the {\it{DSm network architecture}} of the DSm rule of combination.
The first layer of the network consists in all bba of focal elements $A_i, i=1,\ldots,n$ of $m_1(.)$. The second layer of the network consists in all bba of focal elements $B_j, j=1,\ldots,k$ of $m_2(.)$. Each node of layer 2 is connected with each node of layer 1. The output layer (on the right) consists in the combined basic belief assignments of all possible intersections $A_i\cap B_j$, $i=1,\ldots,n $ and $j=1,\ldots,k$. The last step of DSm rule (not included on the figure due to space limitation) consists in the compression of the output layer by regrouping (additioning) all the combined belief assignments corresponding to the same focal elements (by example if $X=A_2\cap B_3=A_4\cap B_5$, then $m(X)=m(A_2\cap B_3)+m(A_4\cap B_5)$). If a third body of evidence provides a new bba $m_3(.)$, the one can combine it by connecting the output layer with the layer associated to $m_3(.)$, and so on. Because of commutativity and associativity properties of DSm rule, the DSm network can be designed with any order for the layers.
The DSm rule of combination can be used for the fusion of any kind of information, whereas the Dempster's rule within Shafer's model can not be used in cases where paradoxist information occurs, or degree of conflict is 1, or when elements of the frame of discernment are not refinable in exclusive finer atoms.

%***************************************
\section{The generation of $D^\Theta$}
%***************************************

\subsection{Memory size requirements}
%****************************************
Before going further on the generation of  $D^\Theta$, it is important to estimate the memory size for storing the elements of $D^\Theta$ for $\vert\Theta\vert=n$. Since each element of $D^\Theta$ can be stored as a $2^n-1$-binary string, the memory size for $D^\Theta$ is given by the right column of the following table (we do not count the size for $\emptyset$ which is 0 and the minimum length is considered here as the byte (8 bits)): 
\begin{center}
\begin{tabular}{l|l|l|l}
$\vert\Theta\vert=n$ & size/elem. & $\#$ of elem. & Size of $\ D^\Theta$\\
\hline
2 &  1 byte &  4 & 4 bytes\\
3 &  1 byte &  18 & 18 bytes\\
4 &  2 bytes &  166 & 0.32 Kb\\
5 &  4 bytes &  7579 & $30$ Kb\\
6 &  8 bytes &  7828352 & $59$ Mb\\
7 &  16 bytes &   $\approx2.4 \cdot {10}^{12}$ & $3.6 \cdot {10}^4$ Gb\\
8 &  32 bytes &  $\approx5.6 \cdot {10}^{22}$ & $1.7 \cdot {10}^{15}$ Gb
\end{tabular}
\end{center}
This table shows the extreme difficulties for our computers to store all the elements of $D^\Theta$ when $\vert \Theta \vert>6$. This complexity must be however relativized with respect to the number of all Boolean functions built from the ultimate refinement (if accessible) $2^{\Theta_{ref}}$ of same initial frame $\Theta$ for applying DST.  The comparison of $\vert D^\Theta \vert$ with respect to
$\vert 2^{\Theta_{ref}} \vert$ is given in the following table
\begin{center}
\begin{tabular}{l|l|l}
$\vert\Theta\vert=n$ & $\vert D^\Theta \vert$ & $\vert 2^{\Theta_{ref}} \vert = 2^{2^n-1}$ \\
\hline
2 &  5 &  $2^3=8$ \\
3 &  19 &  $2^7=128$ \\
4 &  167 &  $2^{15}=32768$ \\
5 &  7580 &  $2^{31}=2147483648$
\end{tabular}
\end{center}
 
\subsection{Monotone Boolean functions}
%*******************************************

 A simple {\it{Boolean function}} $f(.)$ maps $n$-binary inputs $(x_1,\ldots,x_n)\in\{0,1\}^n\triangleq\{0,1\}\times\ldots\times\{0,1\}$ to a single binary output $y=f(x_1,\ldots,x_n)\in\{0,1\}$. Since there are $2^n$ possible input states which can map to either 0 or 1 at the output $y$, the number of possible boolean functions is $2^{2^n}$. Each of these functions can be realized by the logic operations $\wedge$ (and), $\vee$ (or) and $\neg$ (not) \cite{Comtet_1974,Weisstein_2002}. As simple example, let consider only a 2-binary input variable $(x_1,x_2)\in\{0,1\}\times\{0,1\}$ then all the $2^{2^2}=16$ possible Boolean functions $f_i(x_1,x_2)$ built from $(x_1,x_2)$ are summarized in the following tables

\begin{center}
\begin{tabular}{|c|c|c|c|c|c|c|c|c|}
\hline
$(x_1,x_2)$ & $f_0$ & $f_1$ & $f_2$ & $f_3$ & $f_4$ & $f_5$ & $f_6$ & $f_7$\\
\hline
$(0,0)$ & 0 & 0 & 0 & 0 & 0 & 0 & 0 & 0\\
$(0,1)$ & 0 & 0 & 0 & 0 & 1 & 1 & 1 & 1 \\
$(1,0)$ & 0 & 0 & 1 & 1 & 0 & 0 & 1 & 1\\
$(1,1)$ & 0 & 1 & 0 & 1 & 0 & 1 & 0 & 1 \\
\hline
Notation & False & $ x_1\wedge x_2$ & $x_1\wedge \bar{x}_2$ & $x_1$ &  $\bar{x}_1\wedge x_2$ & $x_2$ & $ x_1\veebar x_2$ & $x_1\vee x_2$\\
\hline
\end{tabular}
\end{center}

\begin{center}
\begin{tabular}{|c|c|c|c|c|c|c|c|c|}
\hline
$(x_1,x_2)$ & $f_8$ & $f_9$ & $f_{10}$ & $f_{11}$ & $f_{12}$ & $f_{13}$ & $f_{14}$ & $f_{15}$\\
\hline
$(0,0)$ & 1 & 1 & 1 & 1 & 1 & 1 & 1 & 1 \\
$(0,1)$ & 0 & 0 & 0 & 0 & 1 & 1 & 1 & 1 \\
$(1,0)$ & 0 & 0 & 1 & 1 & 0 & 0 & 1 & 1\\
$(1,1)$ & 0 & 1 & 0 & 1 & 0 & 1 & 0 & 1\\
\hline
Notation &  $x_1 \bar{\vee} x_2$ & $x_1 \triangle x_2$ & $ \bar{x}_2$ & $x_1\vee \bar{x}_2$ &  $\bar{x}_1$ & $\bar{x}_1\vee x_2$ & $ x_1\barwedge x_2$ & True\\
\hline
\end{tabular}
\end{center}

\noindent with the notation $\bar{x}\triangleq\neg x$, $x_1 \veebar x_2\triangleq(x_1\vee x_2)\wedge(\bar{x}_1\vee \bar{x}_2)$ (xor), $x_1 \bar{\vee} x_2\triangleq \neg(x_1\vee x_2)$ (nor), $x_1 \triangle x_2\triangleq (x_1\wedge x_2)\vee(\bar{x}_1\wedge \bar{x}_2)$ (xnor) and
$x_1\barwedge x_2 \triangleq \neg (x_1\wedge x_2)$ (nand). We denote by $\mathcal{F}_n(\wedge,\vee,\neg) =\{f_0(x_1,\ldots,x_n),\ldots, f_{2^{2^n}-1}(x_1,\ldots,x_n)\}$ the set of all possible Boolean functions built from $n$-binary inputs.
Let $\mathbf{x}\triangleq(x_1,\ldots,x_n)$ and $\mathbf{x}'\triangleq({x'}_1,\ldots,{x'}_n)$ be two vectors in $\{0,1\}^n$. Then $\mathbf{x}$ precedes $\mathbf{x}'$ and we denote $\mathbf{x}\preceq\mathbf{x}'$ if and only if $x_i\leq {x'}_i$ for $1\leq i\leq n$ ($\leq$ is applied componentwise). If $x_i < {x'}_i$ for $1\leq i\leq n$ then $\mathbf{x}$ strictly precedes $\mathbf{x}'$ which will be denoted as $\mathbf{x}\prec\mathbf{x}'$. \\

A Boolean function $f$ is said to be a {\it{non-decreasing monotone (or isotone) Boolean function}} (or just {\it{monotone Boolean function}} for short) if and only if $\forall \mathbf{x},\mathbf{x}'\in \{0,1\}^n$ such that $\mathbf{x}\preceq\mathbf{x}'$, then $f(\mathbf{x})\preceq f(\mathbf{x}')$ \cite{Triantaphyllou_2001}.  
Since any isotone Boolean function involves only $\wedge$ and $\vee$ operators (no $\neg$ operations) \cite{Weisstein_2002} and there exists a correspondance between $(\vee,\wedge)$ operators in logics with $(+,\cdot)$ in algebra of numbers and $(\cup,\cap)$ in algebra of sets, the generation of all elements of $D^\Theta$ built from $\Theta$ with $\cup$ and $\cap$ operator is equivalent to the problem of generating isotone Boolean functions over the vertices of the unit $n$-cube. We denote by $\mathcal{M}_n(\wedge,\vee)$ the set of all possible monotone Boolean functions built from $n$-binary inputs. $\mathcal{M}_n(\wedge,\vee)$ is a subset of $\mathcal{F}_n(\wedge,\vee,\neg)$. In the previous example, $f_1(x_1,x_2)$, $f_3(x_1,x_2)$, $f_5(x_1,x_2)$, $f_7(x_1,x_2)$ are isotone Boolean functions but special functions $f_0(x_1,x_2)$ and $f_{2^{2^n}-1}(x_1,\ldots,x_n)$ must also be considered as monotone functions too.  All the other functions belonging to $\mathcal{F}_2(\wedge,\vee,\neg)$ do not belong to $\mathcal{M}_2(\wedge,\vee)$ because they require the $\neg$ operator in their expressions and we can check easily that the monotonicity property $\mathbf{x}\preceq\mathbf{x}' \Rightarrow f(\mathbf{x})\preceq f(\mathbf{x}')$ does not hold for these functions. The Dedekind's problem \cite{Dedekind_1897} is to determine the number $d(n)$ of {\it{distinct}} monotone Boolean functions of $n$-binary variables.  Dedekind \cite{Dedekind_1897} computed $d(0)=2$, $d(1)=3$, $d(2)=6$, $d(3)=20$ and $d(4)=168$. Church \cite{Church_1940} computed $d(5)=7581$ in 1940. Ward \cite{Ward_1946} computed $d(6)=7828354$ in 1946. Church \cite{Church_1965} then computed $d(7)=2414682040998$ in 1965. Between sixties and eighties, important advances have been done to obtain upper and lower bounds for $d(n)$ \cite{Hansel_1966,Kleitman_1969,Korshunov_1981}. In 1991, Wiedemann \cite{Wiedemann_1991} computed $d(8)=56130437228687557907788$ (200 hours of computing time with a Cray-2 processor) which has recently been validated by Fidytek and al. in \cite{Fidytek_2001}. Until now the computation of $d(n)$ for $n>8$ is still a challenge for mathematicians even if the following direct exact explicit formula for $d(n)$ has been obtained by Kisielewicz and Tombak (see \cite{Kisielewicz_1988,Tombak_2001} for proof) :
$$d(n)=\sum_{k=1}^{2^{2^n}} \prod_{j=1}^{2^{n}-1}\prod_{i=0}^{j-1} (1- b_i^k(1-b_j^k)
\prod_{m=0}^{l(i)}(1- b_m^i(1-b_m^j)))$$
\noindent
where $l(0)=0$ and $l(i)=[\log_2 i]$ for $i>0$, $b_i^k\triangleq [k/2^i]-2[k/2^{i+1}]$ and $[x]$ denotes the floor function (i.e. the nearest integer less or equal to $x$). 
The difficulty arises from the huge number of terms involved in the formula, the memory size and the highspeed computation requirements. The last advances and state of art in counting algorithms of Dedekind's numbers can be found in  \cite{Tombak_2001,Fidytek_2001,Triantaphyllou_2001}.

%-------------------------------------------------------------------------
\subsection{Generation of MBF}
%-------------------------------------------------------------------------
Before describing the general algorithm for generating the  monotone Boolean functions (MBF), let examine deeper  the example of section 4.1. From previous tables, one can easily find the set of (restricted) MBF $\mathcal{M}_2^{\star}(\wedge,\vee)=\{f_0(x_1,x_2)=\text{False}, f_1(x_1,x_2)=x_1\wedge x_2, f_5(x_1,x_2)=x_2,f_7(x_1,x_2)=x_1\vee x_2\}$ which is equivalent, using algebra of sets, to hyper-powerset $D^X=\{\emptyset,x_1\cap x_2, x_1,x_2,x_1\cup x_2\}$ associated with frame of discernment $X=\{x_1,x_2\}$. Since the tautology $f_{15}(x_1,x_2)$ is not involved within DSmT, we do not include it as a proper element of $D^X$ and we consider only $\mathcal{M}_2^{\star}(\wedge,\vee)\triangleq \mathcal{M}_2(\wedge,\vee)  \setminus \{f_{15}\}$ rather than $\mathcal{M}_2(\wedge,\vee)$ itself.
Let’s now introduce the Smarandache’s codification for the enumeration of distinct parts of a Venn diagram $X$ with $n$ partially overlapping elements $x_i$,$ i = 1, 2, \ldots, n$.  A such diagram has $2^n-1$ disjoint parts.  One notes with only one digit (or symbol) those parts which belong to only one of the elements $x_i$ (one notes by $<i>$ the part which belongs to $x_i$ only, for $1\leq i\leq n$), with only two digits (or symbols) those parts which belong to exactly two elements (one notes by $<ij>$, with $i<j$, the part which belongs to $x_i$ and $x_j$ only, for $1\leq i <j\leq n$), then with only three digits (or symbols) those parts which belong to exactly three elements (one notes by $<ijk>$ concatenated numbers, with $i<j<k$, the part which belongs to $x_i$, $x_j$, and $x_k$ only, for $1\leq i <j <k \leq n$), and so on up to $<12\ldots n>$ which represents the last part that belongs to all elements $x_i$. For $1\leq n \leq 9$, the Smarandache’s encoding works normally as in base 10.  But, for $n\geq 10$, because there occur two (or more) digits/symbols in notation of the elements starting from 10 on, one considers this codification in base $n+1$, i.e. using one symbol to represent two (or more) digits, for example: $A = 10$, $B = 11$, $C = 12$, etc. For $n =1$ one has only one part, coded $<1>$. For $n =2$ one has three parts, coded $<1>$, $<2>$, $<12>$. Generally, $<ijk>$ does not represent $x_i\cap  x_j\cap x_k$ but only a part of it, the only exception is for $<12\ldots n>$.  For $n =3$ one has $2^3-1= 7$ disjoint parts, coded $<1>$, $<2>$, $<3>$, $<12>$, $<13>$, $<23>$, $<123>$. $<23>$ means the part which belongs to $x_2$ and $x_3$ only, but $<23>\neq x_2\cap x_3$ because $x_2\cap x_3 = \{<23>, <123>\}$ in the Venn diagram of 3 elements $x_1$, $x_2$, and $x_3$.The generalization for $n>3$ is straightforward. Smarandache’s codification can be organized in a numerical increasing order, in lexicographic order or any other orders. An useful order for organizing the Smarandache’s codification for the generation of $D^\Theta$ is the {\it{Dezert-Smarandache order}} $\mathbf{u}_n=[u_1,\ldots,u_{2^n-1}]'$ based on a recursive construction starting with $\mathbf{u}_1\triangleq[<1>]$. Having constructed $\mathbf{u}_{n-1}$, then we can construct $\mathbf{u}_n$ for $n>1$ recursively as follows:
\begin{itemize}
\item include all elements of $\mathbf{u}_{n-1}$ into $\mathbf{u}_{n}$;
\item afterwards, include element $<n>$ as well in $\mathbf{u}_{n}$;
\item then at the end of each element of  $\mathbf{u}_{n-1}$ concatenate the element $<n>$ and
get a new set  ${\mathbf{u}'}_{n-1}$ which then is also included in $\mathbf{u}_{n}$.
\end{itemize}
This is $\mathbf{u}_{n}$, which has $(2^{n-1}-1)+1+(2^{n-1}-1)=2^n-1$ components.
For $n=3$, as example, one gets $\mathbf{u}_3\triangleq[<1>\; <2>\;  <12> \; <3>\; <13> \; <23> \; <123>]'$. Because all elements in $\mathbf{u}_n$ are disjoint, we are able to write each element $d_i$ of $D^X$ in a unique way as a linear combination of $\mathbf{u}_n$ elements, i.e.
\begin{equation}
\mathbf{d}_n=[d_1,\dots,d_{2^n-1}]'=\mathbf{D}_n\cdot \mathbf{u}_n
\end{equation}
Thus $\mathbf{u}_n$ constitutes a basis for generating the elements of $D^X$. Each row in the matrix $\mathbf{D}_n$ represents the coefficients of an element of $D^X$ with respect to the basis $\mathbf{u}_n$.  The rows of $\mathbf{D}_n$ may also be regarded as binary numbers in an increasing order. 
As example, for $n=2$, one has:
\begin{equation}
\underbrace{
\begin{bmatrix}
d_1=x_1\cap x_2\\
d_2=x_2\\
d_3=x_1\\
d_4=x_1\cup x_2
\end{bmatrix}
}_{\mathbf{d}_2}
=
\underbrace{
\begin{bmatrix}
0 & 0 & 1\\
0 & 1 & 1\\
1 & 0 & 1\\
1 & 1 & 1\\
\end{bmatrix}
}_{\mathbf{D}_2}
\cdot
\underbrace{
\begin{bmatrix}
<1>\\
<2>\\
<12>
\end{bmatrix}
}_{\mathbf{u}_2}
\label{eq:D2}
\end{equation}

\noindent where the "matrix product" is done after identifying $(+,\cdot)$ with $(\cup,\cap)$, $0\cdot <x>$ with $\emptyset$ and $1\cdot <x>$ with $<x>$. The generation of $D^X$ is then strictly equivalent to generate $\mathbf{u}_n$ and matrix $\mathbf{D}_n$ which can be easily obtained by the following recursive procedure:
\begin{itemize}
\item start with  $\mathbf{D}_0^c=[0 \, 1]'$  corresponding to all Boolean functions with no input variable ($n=0$).
\item build the $\mathbf{D}_1^c$ matrix from each row $\mathbf{r}_i$ of $\mathbf{D}_0^c$ by adjoining it  to any other row $\mathbf{r}_j$ of $\mathbf{D}_0^c$ such that $\mathbf{r}_i \cup \mathbf{r}_j=\mathbf{r}_j$. This is equivalent here to add either 0 or 1 in front (i.e. left side) of $\mathbf{r}_1\equiv 0$ but only 1 in front of $\mathbf{r}_2\equiv 1$.  Since the tautology is not involved in hyper-powerset, then one has to remove the first column and the last line of 
\begin{equation*}
\mathbf{D}_1^c=\begin{bmatrix}
0 & 0\\
0 & 1\\
1 & 1\\
\end{bmatrix}
\quad \text{to obtain finally} \quad 
\mathbf{D}_1=\begin{bmatrix}
0\\
1\\
\end{bmatrix}
\end{equation*}
\item build $\mathbf{D}_2^c$ from $\mathbf{D}_1^c$ by adjoining to each row $\mathbf{r}_i$ of $\mathbf{D}_1^c$, any row $\mathbf{r}_j$ of $\mathbf{D}_1^c$ such that $\mathbf{r}_i \cup \mathbf{r}_j=\mathbf{r}_j$ and then remove the first column and the last line of $\mathbf{D}_2^c$ to get $\mathbf{D}_2$ as in \eqref{eq:D2}.
\item build $\mathbf{D}_3^c$ from $\mathbf{D}_2^c$ by adjoining to each row $\mathbf{r}_i$ of $\mathbf{D}_2^c$ any row $\mathbf{r}_j$ of $\mathbf{D}_2^c$ such that $\mathbf{r}_i \cup \mathbf{r}_j=\mathbf{r}_j$ and then remove the first column and the last line of $\mathbf{D}_3^c$ to get $\mathbf{D}_3$ given by
\begin{equation*}
\mathbf{D}_3=\begin{bmatrix}
    0  &   0    & 0   &  0   &  0   &  0   &  0\\
     0 &    0   &  0  &   0  &   0  &   0  &   1\\
     0  &   0   &  0   &  0  &   0  &   1   &  1\\
     0  &   0  &   0  &   0  &   1  &   0   &  1\\
     0  &   0  &   0  &   0  &   1  &   1   &  1\\
     0   &  0  &   0  &   1  &   1  &   1   &  1\\
     0  &   0   &  1   &  0  &   0  &   0   &  1\\
     0  &   0  &   1  &   0  &   0  &   1   &  1\\
     0  &   0  &   1  &   0  &   1  &   0  &   1\\
     0  &   0   &  1  &   0  &   1  &   1  &   1\\
     0  &   0   &  1  &   1   &  1  &   1  &   1\\
     0  &   1  &   1 &    0  &   0   &  1  &   1\\
     0  &   1  &   1  &   0  &   1   &  1  &   1\\
     0  &   1  &   1  &   1   &  1  &   1  &   1\\
     1  &   0  &   1  &   0  &   1  &   0 &    1\\
     1  &   0  &   1  &   0  &   1  &   1   &  1\\
     1  &   0  &   1  &   1  &   1  &   1 &    1\\
     1   &  1  &   1  &   0   &  1  &   1  &   1\\
     1  &   1  &   1  &   1  &   1  &   1  &   1
\end{bmatrix}
\end{equation*}
\item Likewise, $\mathbf{D}_n^c$ is built from $\mathbf{D}_{n-1}^c$ by adjoining to each row $\mathbf{r}_i$ of $\mathbf{D}_{n-1}^c$ any row $\mathbf{r}_j$ of $\mathbf{D}_{n-1}^c$ such that $\mathbf{r}_i \cup \mathbf{r}_j=\mathbf{r}_j$. Then $\mathbf{D}_n$ is obtained by removing the first column and the last line of $\mathbf{D}_n^c$.
\end{itemize}

For convenience, we provide here the source code in Matlab language to generate $D^\Theta$.
This code includes the identification of elements of $D^\Theta$ corresponding to each monotone Boolean  function according to the Smarandache's codification. 
\begin{center}
\lstset{language=Matlab,stepnumber=1,showspaces=false,numbersep=5pt,numberstyle=\tiny,basicstyle=\scriptsize}
\begin{lstlisting}
%************************************************
% Copyright (c)  2003 J.Dezert and F.Smarandache
% 
% Purpose: Generation of D^Theta for the DSmT for  
% Theta={theta_1,..,Theta_n}. Due to the huge 
% # of elements of D^Theta only cases up to n<7  
% are usually tractable on computers.     
%************************************************
n=input('Enter cardinality for Theta (0<n<6) ?');
% Generation of the Smarandache codification
% Note: this should be implemented using
% character strings for n>9
u_n=[1];
for nn=2:n
u_n=[u_n nn (u_n*10+nn*ones(1,size(u_n*10,2)))];
end
% Generation of D_n (isotone boolean functions)
D_n1=[0;1];
for nn=1:n, D_n=[];
    for i=1:size(D_n1,1),Li=D_n1(i,:);
	for j=i:size(D_n1,1)
	 Lj=D_n1(j,:);Li_inter_Lj=and(Li,Lj);
	 Li_union_Lj=or(Li,Lj);
	 if((Li_inter_Lj==Li)&(Li_union_Lj==Lj))
		D_n=[D_n; Li Lj];
	 end
	end
    end
    D_n1=D_n;
end
DD=D_n;DD(:,1)=[];DD(size(DD,1),:)=[]; D_n=DD;
% Result display
disp(['|Theta|=n=',num2str(n)])
disp(['|D^Theta|=',num2str(size(D_n,1))])
disp('Elem. of D^Theta are obtained by D_n*u_n')
disp(['with u_n=[',num2str(u_n),']'' and'])
D_n=D_n

\end{lstlisting}
\vspace{-0.3cm}
\begin{center}
\bf{Matlab\footnote{Matlab is a trademark of The MathWorks, Inc.} source code for generating $D^\Theta$}
\end{center}
\end{center}

\section{Conclusion}

The DSmT proposes a new solution to combine conflicting sources of information in some problems  where the frame of discernment $\Theta=\{\theta_1,\ldots,\theta_n\}$ cannot be considered as an exhaustive and exclusive finite set of hypotheses. The DSmT deals with elements $\theta_i$ which have possibly (but not necessarily) continuous and/or relative interpretation to  the corpus of evidences (like the notions of smallness/tallness, beauty/ugliness, pleasure/pain, heat/coldness, even the notion of colors - due to the continuous spectrum of the light, $\ldots$); the interpretation of $\theta_i$ through the bba mechanism given by each source being, in general, built only from its own limited knowledge/experience and senses. 
This DSm model can be considered as the opposite of the Shafer's model on which is based the DST. The DSmT, based on the notion of hyper-powerset $D^\Theta$
over $\Theta$ and the refutation of the third middle excluded, requires in theory to manipulate the basic beliefs assigned of every element of $D^\Theta$. A powerful method and a source code to generate recursively all the elements of $D^\Theta$ has been presented in this paper to help the reader to solve a wide class of fusion problems with the DSmT.

\end{document}